\documentclass{amsart}
\usepackage{ifthen, url}

\newcommand{\from}{\colon}
\newcommand{\Nat}{\mathbb{N}}

\newtheorem{axiom}{Axiom}
\newtheorem{theorem}{Theorem}
\newtheorem{lemma}{Lemma}
\newtheorem{proposition}{Proposition}
\newtheorem{corollary}{Corollary}
\theoremstyle{definition}
\newtheorem{definition}{Definition}
\newtheorem{example}{Example}

\theoremstyle{remark}

\newcommand{\theor}[3][]{\begin{theorem}[#1] #3 \label{thm:#2} \end{theorem}}
\def\thm#1{{t}heorem \ref{thm:#1}}

\def\lema#1#2{\begin{lemma} #2 \label{lem:#1} \end{lemma}}
\def\lem#1{{l}emma \ref{lem:#1}}

\def\display#1#2{\begin{equation} #2 \label{eqn:#1} \end{equation}}
\def\eqn#1{(\ref{eqn:#1})}

\newboolean{RANW}
\setboolean{RANW}{false}
\newlength{\testWidth} \newlength{\badWidth}

\newcommand{\useRsltAndNumberwithin}{
\setboolean{RANW}{true}
\settowidth\badWidth{\ref{useRsltAndNumberwihin_generates_spurious_undefined-reference_messages.}}
\newcommand{\badMatch}{\equal{\the\testWidth}{\the\badWidth}}
}

\newcommand{\rslt}[1]{%
\ifthenelse{\boolean{RANW}}%
{%
{\settowidth{\testWidth}{\ref{lem:#1}}\ifthenelse{\badMatch}%
{%
\settowidth{\testWidth}{\ref{prp:#1}}\ifthenelse{\badMatch}%
{%
\settowidth{\testWidth}{\ref{thm:#1}}\ifthenelse{\badMatch}%
{%
\settowidth{\testWidth}{\ref{cor:#1}}\ifthenelse{\badMatch}{\textbf{result ??}}%
{corollary \ref{cor:#1}}%
}%
{theorem \ref{thm:#1}}%
}%
{proposition \ref{prp:#1}}%
}%
{lemma \ref{lem:#1}}%
}}%
{%
\ifthenelse{\ref{lem:#1} > 0}{lemma \ref{lem:#1}}%
{%
\ifthenelse{\ref{prp:#1} > 0}{proposition \ref{prp:#1}}%
{%
\ifthenelse{\ref{thm:#1} > 0}{theorem \ref{thm:#1}}%
{%
\ifthenelse{\ref{cor:#1} > 0}{corollary \ref{cor:#1}}{\textbf{Result ??}}%
}%
}%
}%
}}

\newcommand{\eqnref}[1]{\ref{eqn:#1}}
 
\newcommand{\rsltref}[1]{%
\ifthenelse{\boolean{RANW}}%
{%
{\settowidth{\testWidth}{\ref{lem:#1}}\ifthenelse{\badMatch}%
{%
\settowidth{\testWidth}{\ref{prp:#1}}\ifthenelse{\badMatch}%
{%
\settowidth{\testWidth}{\ref{thm:#1}}\ifthenelse{\badMatch}%
{%
\settowidth{\testWidth}{\ref{cor:#1}}\ifthenelse{\badMatch}{\textbf{result ??}}%
{\ref{cor:#1}}%
}%
{\ref{thm:#1}}%
}%
{\ref{prp:#1}}%
}%
{\ref{lem:#1}}%
}}%
{%
\ifthenelse{\ref{lem:#1} > 0}{\ref{lem:#1}}%
{%
\ifthenelse{\ref{prp:#1} > 0}{\ref{prp:#1}}%
{%
\ifthenelse{\ref{thm:#1} > 0}{\ref{thm:#1}}%
{%
\ifthenelse{\ref{cor:#1} > 0}{\ref{cor:#1}}{\textbf{Result ??}}%
}%
}%
}%
}}

\newenvironment{enum} 
{\begin{list}{\makebox[\labelwidth][l]{(\arabic{enumi})}}{\usecounter{enumi}}
\setcounter{enumi}{\value{equation}}}
{\setcounter{equation}{\value{enumi}} \end{list}}

\newcommand{\meti}[2]{\item #2 \label{eqn:#1}} 

\newcommand{\authorGreen}{\ifthenelse{\isundefined{\authorInPageHeader}}
{\author[\relax]{Edward J. Green}}
{\author[\authorInPageHeader]{Edward J. Green}}}

\authorGreen
\address{Department of Economics, The Pennsylvania State University,
  University Park, PA 16802, USA}
\email{eug2@psu.edu}

\renewcommand{\b}[1]{\mathbf{#1}}

\title[Individual-Level Randomness]{Individual-Level Randomness\\ in a Nonatomic Population}

\thanks{The research reported here was initiated while I was
  affiliated with the University of Minnesota and the Federal Reserve
  Bank of Minneapolis. A significant part of it was conducted during a
  visit to the Mathematical Sciences Research Institute in December,
  1985.  I would like to express my appreciation for the support and
  hospitality of the Institute.  I would particularly like to thank
  William Zame for his comments and suggestions, and Jerome Keisler
  and Robert Anderson for their correspondence regarding the
  alternative approach, via Loeb measure, to this problem.}

\date{April, 2019}



\begin{document}

\begin{abstract}
This paper provides a construction of an uncountable family of
i.i.d.~random vectors, indexed by the points of a nonatomic measure
space, such that (a) a sample is a measurable function from the index
space, and (b) an idealization of the Glivenko-Cantelli theorem holds
exactly with respect to the measure on that space.  That is, samples
possess a.s.~the distribution from which they are drawn.  Moreover,
any subspace of the index space with positive measure inherits the
same property.  This homogeneity property is important for
applications of the construction in economics.
\end{abstract}

\maketitle

\section{Introduction}

The Glivenko-Cantelli theorem states that, almost surely, the sample
distributions of i.i.d.~random variables converge weakly to the
statistical distribution of the random variables.  (Parthasar\-athy
(1967), Theorem II.7.1.)  This paper provides a set-theoretic
construction of an uncountable family of i.i.d.~random vectors,
indexed by the points of a nonatomic measure space, such that (a)
samples are measurable functions from the index space, and (b) an
exact analogue of the Glivenko-Cantelli theorem holds with respect to
the measure on that space.  That is, a sample can be viewed as a
random vector by regarding the index space itself as a probability
space, and a.s.~the sample possesses the same distribution as that of
the i.i.d.~random vectors from which it is drawn.  Moreover, any
subspace of the index space with positive measure inherits the same
property, if the measure of the subspace is normalized to be a
probability measure.  This homogeneity property is important for an
application of the construction in economics which will
be discussed below. The construction presented here is an alternative
to the construction via Loeb measure, first presented by Keisler
(1977) and subsequently simplified by Anderson (1991).

To understand what this construction accomplishes, first consider
a more direct construction of a family of i.i.d.~random variables
for which an exact idealization of the strong law of large numbers
holds.
This construction begins with Kolmogorov's construction of a
continuum of i.i.d.~random variables
$\{ \phi_{t}  \mid  t\in [0,1]\}$ having a prescribed
distribution with finite first moment.
Then the measure on the sample space $\Omega$ so constructed
is extended in such a way that
$\int_{0}^{1}\phi_{t}(\omega )dt=E(\phi_{0})$ a.s. 
This equation idealizes the strong law of large numbers, with
Lebesgue integration on $[0,1]$ playing the role of
averaging over the sample.
The techniques needed to prove the existence of a continuous-time
i.i.d.~process satisfying this integral equation were developed
by Doob (cf. (1937), (1947), (1953) chapter II), and an existence
proof has been given in full by Judd (1985).

This construction has been widely cited by economic theorists.  In
particular, the i.i.d.~processes just discussed have been thought to
provide the mathematical basis for tractable models of economies in
which individual traders face idiosyncratic risks---risks of gains or
losses that are sizable for each individual trader, but that are
independent across traders and accurately predictable in the
aggregate.  (The risk of death is tolerably close to being
idiosyncratic, at least in populations where epidemic diseases are
under control, but the risk of outbreak of war is not idiosyncratic.
This contrast explains why life insurance is easily available but
insurance against political risks is not.)  In such models, $[0,1]$
represents the population of traders, $\phi_{t}$ represents the random
gain or loss experienced by trader $t,$ and
$\int_{0}^{1}\phi_{t}(\omega )dt$ represents the aggregate net gain or
loss to the economy {\it ex post} when $\omega$ is the state of the
world.

Feldman and Gilles (1985) provide documentation of the
importance of models of this genre to current economic theory.
They also argue that the construction just described is actually
inadequate to provide a mathematical foundation for those models.
The problem is that the models posit more than just the one
integral equation $\int_{0}^{1}\phi_{t}(\omega)dt =
E(\phi_{0})$ a.s. 
It is assumed that, for every $\theta\in
(0,1],\int_{0}^{\theta}\phi_{t}(\omega)dt=\theta E(\phi_{0})$.
This homogeneity assumption reflects the economic idea that any
non-negligible fraction of the traders in a large economy could
potentially form a risk-pooling coalition that would provide its
members with virtually complete insurance against idiosyncratic
risks, and that therefore the per-capita aggregate resources {\it
ex post} of a large coalition in an economy of traders facing 
i.i.d.~risks should not depend on the coalition.
In the study of insurance, this assumption is required in order
to demonstrate that equal sharing of resources is the unique
cooperative arrangement (technically, the unique core allocation)
in this economy.

Mathematically this strengthened assumption may seem innocuous,
because the interval $[0,\theta ]$ with normalized
Lebesgue measure is isomorphic to $[0,1]$.
Thus the integral equation should be as plausible for any value
of $\theta$ as it is for $\theta = 1$.
However, Feldman and Gilles have shown that it is inconsistent to
make the assumption for all $\theta ,$ if $\phi$ is a
Bernoulli process taking values {\it always} in 
$\{ 0,1\}$.\footnote{There are applications in which it is
crucial that $\phi$ should take values only in
$\{0,1\}$.
For instance, if $\phi$ is the characteristic function of
some event, then 0 and 1 are the only values that it can
meaningfully take.}
In that case, the assumption is equivalent to
$\int_{0}^{\theta}\phi_{t}(\omega )dt=\theta / 2$ a.s.
Considering this equation for rational values of $\theta$ (to
assure measurability of the event that the equation holds for all
values considered), the Radon-Nikodym theorem implies that a.s.
$\{t  \mid  \phi_{t}(\omega ) = 1/2\}$ has Lebesgue measure 1.
This conclusion contradicts the restriction that has been imposed
on the range of $\phi ,$ though.

In view of this contradiction, the body of economic theory that
has been formulated in terms of this model needs to be placed on
a more secure foundation.\footnote{Feldman and Gilles (1985) cite a
number of prominent contributions to economic theory that their
argument shows to be inconsistent (at least if risk is
parametrized as being Bernoulli).
These authors, as well as Judd (1985), Uhlig (1996), and others,
have proposed alternative definitions of integrals over the
population in order to assure homogeneity.
However, besides the general disadvantages of recourse to such
alternatives cited by Doob (1947), there are economic arguments
that employ Lebesgue integration with respect to countably
additive measures on both the sample space and the population.
(An example is Green (1987), where the results proved here are
used.)
This requirement motivates the present study.}
This task will be accomplished here.
It will be proved that there exists an i.i.d.~family of random
vectors that satisfy an {\it analogue} of the condition that, for
each measurable set $A$ in the range of $\phi_{0}$ and
for each $\theta\in [0,1],$
$\lambda (\{ t \mid \phi_{t}(\omega )\in A\}\bigcap[0,\theta ])
= \theta P(\phi_{0}\in A)$ a.s.
(Here $\lambda$ denotes Lebesgue measure on $[0,1]$ and
$P$ denotes probability measure defined on the sample space.)
This condition just stated would be an exact idealization of the
Glivenko-Cantelli theorem,
but in the case of integrable, $\{0,1\}$-valued random
variables it would imply the set of integral equations that has
just been seen to be inconsistent.
This problem will be avoided by indexing random vectors by
elements of an abstract nonatomic probability space rather than
by numbers in the unit interval.
Correspondingly, sample functions will be integrated over
measurable subsets of this space rather than over intervals.
There is no distinction of economic realism between numbers in
the unit interval and sample points of another probability space
as names of idealized traders.
Thus, the generalized stochastic process described here provides
just as appropriate an economic model as does a process indexed
by the unit interval.

The key to avoiding the contradiction derived by Feldman and
Gilles is to represent the population as a probability space
having a $\sigma$-algebra that is not countably generated.
Then the strategy of restricting attention to the countable set
of intervals $\{[0,\theta ] \mid \theta \text{ is rational} \}$ and
subsequently appealing to the Radon-Nikodym theorem cannot be
emulated.
Rather than starting with a given sample space $\Omega$ and
population space $\Theta ,$ a function $\phi$ will be
defined on the Cartesian product of arbitrary sets $\Omega$
and $\Theta ,$ and then these sets will be endowed with
probability structure in a way that guarantees the required
properties of $\phi $.
Consequently the function $\phi$ and the sets $\Omega$
and $\Theta$ with their respective $\sigma$-algebras will
constitute a universal limit process.
When appropriate probability measures are defined on these
$\sigma$-algebras, the sections of $\phi$ will become an 
i.i.d.~family of random vectors having any specified distribution,
and satisfying the idealized Glivenko-Cantelli property relative
to any measurable subset of the population.

Before carrying out the details of this program, two features of
the limit process need to be discussed.
First, it would be mathematically possible to endow the sample
space and the population space with any combination of
probability measures.
This means that it would be possible for all of the random
vectors to have one distribution, but for the sample functions 
a.s.~to have another distribution.
There seems to be no way to rule out this possibility by appeal
to the kinds of separability or joint-measurability
considerations that are usually invoked.
However, the interpretation of the process as a limit object
provides a strong reason to impose a connection between the two
measures.
This is made clear in the next section.

Second, although all of the sections of the function $\phi$
(defined by fixing either a sample point or a member of the
population) are measurable, $\phi$ is not jointly measurable
in its two arguments.
At the end of the paper, it will be proved that no limit process
for a nondegenerate distribution can be jointly measurable.
Joint measurability might have been of interest for two reasons.
First, it could serve as a selection criterion that would
eliminate counterintuitive processes from consideration.
It has been mentioned above that these counterintuitive processes
can be identified as being pathological on other grounds.
Second, joint measurability would justify the application of the
Fubini theorem to the process.
However, the marginal distributions of the process are a.s.
constant, and the constants for the two variables are equal, so
the conclusion of the Fubini theorem is satisfied even though the
process is not jointly measurable.
Thus the failure of joint measurability is not a serious problem
here.

\section{The sample-distribution limit of an i.i.d.~sequence}

Suppose that $\b\Omega = (\Omega ,\mathcal{B},\pi )$ and
$\b{R}=(R,\mathcal{R},\mu )$ are probability spaces and that
$\{\phi_{n} \from \Omega\rightarrow R\}_{n\in\Nat}$ is a sequence of
independent random vectors, each having distribution $\mu $.  The
finite sample $(\phi_{k}(\omega ))_{k<n}$ can be regarded as a random
vector on the probability space $\Theta_{n}=(\{0,\ldots
,n-1\},\mathcal{F}_{n},\nu_{n}),$ where $\mathcal{F}_{n}$ is the power
set of $\{0,\ldots ,n-1\}$ and $\nu_{n}$ is normalized counting
measure.  In the case that $\b{R}$ is a separable metric space, the finite
sample distributions $\nu_{n}$ converge weakly to $\mu$ almost surely
as $n$ tends to infinity.  (Parthasarathy (1967), Theorem II.7.1.)  In
view of this fact, it is natural to look for a probability space
$\b\Theta =(\Theta,\mathcal{F},\nu )$ and a function $\phi  \from \Omega\times
\Theta\rightarrow R$ such that\footnote{$\phi_{\omega}(\theta
  )=\phi_{\theta}(\omega )=\phi(\omega ,\theta )$.}

\begin{enum}
\meti a{$\{\phi_{\theta} \from \Omega\rightarrow
R\}_{\theta\in\Theta}$ are independent random vectors
having distribution $\mu$, and}
\meti  b{the sample functions 
$\phi_{\omega} \from \Theta\rightarrow R$ are measurable and have
  distribution $\mu$ a.s.}
\end{enum}

An infinite family of random variables is independent iff every
finite subfamily is independent.
``Almost surely'' will always refer to events in $\mathcal{B}$
rather than in $\mathcal{F}$.
 If $\b\Phi =(\b\Omega ,\b\Theta ,\phi)$ satisfies
these conditions, then it will be called a 
\emph{sample-distribution limit} for $\b{R}$.
This paper proves the existence of a sample-distribution limit
that is homogeneous in the sense that restricting $\nu$ to
any set $A\in \mathcal{F}$ such that
$\nu (A) > 0$ and normalizing so that
$\nu (A)=1$ yields again a sample-distribution limit for
$\b{R}$.
(The impossibility of this, if $\b\Theta$ is the unit interval
with Lebesgue measure, is the result of Feldman and Gilles (1985)
that has been discussed above.)
However it will be proved here that, subject to a mild
restriction, no sample-distribution limit is measurable with
respect to $\mathcal{B}\times\mathcal{F}$.

The construction of a homogeneous sample-distribution limit will
rely heavily on Kolmogorov's construction of a set of independent
random vectors having distribution
$\b{R}=(R,\mathcal{R},\mu )$.
The relevant details of Kolmogorov's construction are now
reviewed.
Throughout this paper it will be assumed that

\display c{\mathcal{R}\neq \{\emptyset ,R\} \text{,  and }
\Theta \text{ and }\Omega \text{ are disjoint infinite sets.}}

Kolmogorov's construction takes $R^{\Theta}$ to be the sample space.
To define a $\sigma$-algebra and probability measure, first define

\display d{\mathcal{A}_{\mathcal{B}}^{0}=\{ \alpha \mid \alpha
  \from \Theta \to \mathcal{R}\setminus \{\emptyset\}
  \text{ and } \Theta \setminus \alpha^{-1}(R) \text{ is finite}\}.}

\noindent Every $\alpha\in\mathcal{A}_{\mathcal{B}}^{0}$ can be
regarded as specifying a subset $\alpha^{\#}$ of
$R^{\Theta}$ by

\display e{x\in \alpha^{\#} \text{ iff } \forall\theta \; x(\theta
) \in \alpha(\theta ).}

\noindent The sets $\alpha^{\#}$ defined by \eqn{e} are called
{\it cylinders.}
For $X\subseteq R^{\Theta},$ define
$A\subseteq\mathcal{A}_{\mathcal{B}}^{0}$ to be a {\it cylindrical
partition} of $X$ iff

\display f{A \text{ is a finite subset of } \mathcal{A}_{\mathcal
B}^{0} \text{ and } \{\alpha^{\#}\}_{\alpha \in A} \text{ is a
partition of } X.}

\noindent Define $\mathcal{A}_{\mathcal{B}}$ to be the set of subsets
of $R^{\Theta}$ having a cylindrical partition.
That is,

\display g{X\in \mathcal{A}_{\mathcal{B}} \text{ iff } \exists A \; [A
\text{ satisfies \eqn{f} w.r.t.}~X]}.

\noindent Now, begin to define the Kolmogorov extension measure
$\kappa_{\pi},$ by defining it on $\mathcal{A}_{\mathcal{B}}$.
\medskip
\display h{\kappa_{\pi}(X)=\Sigma_{\alpha \in A}[\Pi_{\theta\in\Theta}\mu
(\alpha (\theta ))] \text{ if $A$ satisfies \eqn{f} w.r.t.~$X$.}}

\noindent Define $\mathcal{E}_{\mathcal{B}}$  by

\display C{\mathcal{E}_{\mathcal{B}} \text{ is the smallest
    $\sigma$-algebra containing } \mathcal{A}_{\mathcal{B}}.}

\noindent The definition of $\kappa_{\pi}$ will be extended to
$\mathcal{E}_{\mathcal{B}}$.
The following lemma summarizes results of a series of arguments
and constructions that are described in Halmos (1974), \S 33, 
\S 37, \S 38.\footnote{Halmos deals explicitly only with the case that
$\Theta =\Nat,$ but the argument is completely general.}

\begin{lemma}[Kolmogorov]  The definition of
$\kappa_{\pi}(X)$ in \eqn{h} does not depend on which
cylindrical partition of $X$ is used.
$\mathcal{A}_{\mathcal{B}}$ is an algebra of subsets of 
$R^{\Theta}$.
There is a unique probability measure $\kappa_{\pi}$ defined
on $\mathcal{E}_{\mathcal{B}}$ that satisfies \eqn{h} for every 
$X\in \mathcal{A}_{\mathcal{B}}$ and for every $A$ that satisfies
\eqn{f} w.r.t.~X. 

The projections $p_{\theta} \from R^{\Theta}\rightarrow R$ defined
by
\display i{p_{\theta}(x) = x(\theta) \text{ are independent random
    vectors having distribution } \mu.}
\label{lem:a} \end{lemma}

Note that the last assertion follows directly from \eqn{h}.

It is useful to know that a set in $\mathcal{E}_{\mathcal{B}}$ is
defined in terms of restrictions on only a countable set of
coordinates.
That is,

\lema b{If $X\in \mathcal{E}_{\mathcal{B}},$ then
there exists $\Theta_{X}\subseteq \Theta$ such that
\display j{\begin{aligned} \Theta_{X} \text{ is countable and }&\\
  \forall y\in R^{\Theta} \;
[y\in X\iff &[\exists x\in X \; \forall \theta\in \Theta_{X} \; [x(\Theta )
= y(\Theta )]]. \end{aligned}}}

\begin{proof} It is easily verified that the set of all subsets of
$R^{\Theta}$ for which there exists 
$\Theta_{X}\subseteq \Theta$ satisfying \eqn{j} is a
$\sigma$-algebra containing $\mathcal{A}_{\mathcal{B}}$.
Thus $\mathcal{E}_{\mathcal{B}}$ is a sub $\sigma$-algebra, since
it is the smallest $\sigma$-algebra containing
$\mathcal{A}_{\mathcal{B}}$. \end{proof}

\section{Construction of a sample-distribution limit from a rich
function}

Define a \emph{rich function} to be a function
$\phi  \from \Omega\times\Theta\rightarrow R$ that satisfies 

\display k {\begin{split}
  \forall h\in R^\Nat\;
[\forall f\in\Theta^\Nat\exists \omega \; [f \text{ is 1--1 } \implies  \forall
n\phi (\omega,f(n)) = h(n)] \text{ and }&\\
\forall g\in\Omega^\Nat\exists \theta \;[g \text{ is 
    1--1 } \implies  \forall n\phi (g(n),\theta ) = h(n)]].&
\end{split}}

\noindent Note that $f,\Theta$ and $\omega$ are dual to 
$g$, $\Omega$ and $\theta$ in \eqn{k}.

Given a rich function $\phi$ and a probability space $\b{R}$, $\b\Omega$
is now constructed.  The idea that guides this construction is that
\eqn{k} guarantees enough diversity in the behavior of $\phi$ on
$\Omega$ so that a $\sigma$-algebra isomorphic to that constructed by
Kolmogorov is required to make every $\phi_{\theta}$ measurable.  This
isomorphism commutes in an appropriate sense with $\phi$ and the
corresponding stochastic process \eqn{i} of Kolmogorov's construction.
Using this fact, Kolmogorov's construction can be pulled back to
$\Omega$ to define $\b\Omega$ in such a way that the $\phi_{\theta}$
are i.i.d.~with distribution $\mu $.  Then, by the duality in \eqn{k},
$\b\Theta$ can be constructed analogously so that the $\phi_{\omega}$
are i.i.d.~with distribution $\mu $.  Thus $\b\Phi$ is a
sample-distribution limit for $\b{R}$.  The independence of the random
vectors $\phi_{\omega}$ will be further exploited to show that
$\b\Phi$ is homogeneous.

$\Omega$ is given as a Cartesian factor space of the domain of
the rich function $\phi $.
$\mathcal{B}$ is now defined as the range of a mapping
$\Psi_\Omega  \from \mathcal{E}_{\mathcal{B}}\rightarrow \mathcal{P}(\Omega ),$ where
$\mathcal{E}_{\mathcal{B}}$ continues to denote the
$\sigma$-algebra on $R^{\Theta}$ obtained in \lem{a}.
Specifically,
\display l{\forall X\in \mathcal{E}_{\mathcal{B}} \enspace \Psi_\Omega (X) = \{ \omega \mid \phi_{\omega}\in X\}}
\noindent and
\display m{\mathcal{B}=\Psi_\Omega (\mathcal{E}_{\mathcal{B}}).}

\lema c{If $\phi$ is rich, then ${\mathcal
B}$ is a $\sigma$-algebra of subsets of $\Omega ,$ and 
$\Psi_\Omega  \from \mathcal{E}_{\mathcal{B}}\rightarrow \mathcal{B}$ is an
isomorphism of $\sigma$-algebras.}

\begin{proof} It is clear that $\mathcal{B}$ is a $\sigma$-algebra and
that $\Psi_\Omega$ is a homomorphism onto $\mathcal{B}$.
To show that $\Psi_\Omega$ is an isomorphism, it is only necessary
to show that it is $1-1$.
Suppose that $X\in\mathcal{E}_{\mathcal{B}}$ and $Y\in {\mathcal
E}_{\mathcal{B}}$ and $X\neq Y$.
Without loss of generality, assume that
$x\in X\setminus Y$.
Let $\Theta_{X}$ and $\Theta_{Y}$ be countable sets
possessing the property \eqn{j} for $X$ and $Y$
respectively, which are guaranteed to exist by \lem{b}, and let
$f \from \Nat\rightarrow\Theta$ be $1-1$ and
$\Theta_{X}\cup\Theta_{Y}\subseteq f(\Nat)$.
Define $h \from \Nat\rightarrow R$ by 
$h(n)=x(f(n))$.
Then, by \eqn{k}, there exists an $\omega$ such that
$\forall n\;\phi (\omega ,f(n)) =h(n)$.
By \eqn{j}, then, $\phi_{\omega}\in X\setminus Y$.
By \eqn{l}, $\omega\in\Psi_\Omega (X)\setminus\Psi_\Omega (Y),$ so 
$\Psi_\Omega (X)\neq \Psi_\Omega (Y)$. \end{proof}

In view of \lem{c}, it is clear how to define $\pi $.
Namely,

\display n{\pi (X) = \kappa_{\pi} (\Psi_\Omega^{-1}(X)).}

\lema d{If $\phi$ is rich, then
$\b\Omega = (\Omega ,\mathcal{B} ,\pi )$ defined by 
\eqn{l} - \eqn{n} is a probability space on which the random
vectors $\phi _{\theta}$ are i.i.d.~with distribution $\b{R}$.}

\begin{proof} \uppercase\lem{a} and \lem{c} show that $\b\Omega$ is a
probability space.  The independence assertion also follows from
\lem{a}, using the equations $\phi_{\theta}(\omega ) =
p_{\theta}(\phi_{\omega})$ and $\Psi_\Omega (\alpha^{\#})=\{\omega \mid \forall
\theta\phi_{\omega}(\theta )\in\alpha (\theta )\} ,$ which are
consequences of the definitions above. \end{proof}

As mentioned after equation \eqn{k}, this argument can be dualized
with respect to $\b\Omega$ and $\b\Theta $.
Let $(n')$ be the dual of equation $(n),$ and let Lemma
$n'$ be the dual of Lemma $n$.

\noindent The following result is an immediate consequence of
\lem{d} and lemma \rsltref{d}$'$.

\theor A{If $\phi$ is rich, then
$\b\Phi =(\b\Omega,\b\Theta ,\phi )$ is a sample-distribution
limit for $\b{R}$.}

\section{Homogeneity of $\phi$}

Define $\b\Phi$ to be a {\it homogeneous} sample-distribution
limit for $\b{R}$ if it is a sample-distribution limit for that
distribution and also

\display o{\forall A\in\mathcal{F}\quad \forall B\in{\mathcal
R}\quad [\nu (A\cap\phi_{\omega}^{-1}(B)) = \nu (A)\mu(B)
\text{ a.s.}]}

\noindent If positive-measure sets in $\mathcal{F}$ are
analogous to infinite subsets of $\Nat,$ then \eqn{o}
intuitively ought to hold in the limit because the
sequential-convergence result cited at the beginning of the paper
applies to every infinite subsequence of
$\{\phi_{n}\}_{n\in\Nat}$.

It will now be shown that the sample-distribution limit
$\b\Phi$ just constructed is homogeneous, that is, that
$\phi$ satisfies \eqn{o}.
For any $A\in\mathcal{F}$ and $B\in \mathcal{R},$ there will 
be a countable subset of $\Omega$ where the condition
asserted by \eqn{o} to hold a.s.~is violated.
Thus, it must be shown that a countable subset is a subset of a
probability-zero event of $\mathcal{B}$.
The next two lemmas establish this.

\lema e{If $\phi$ is rich, $r\in R,$ and
$\omega\in \Omega ,$ then the cardinal of 
$\{\theta \mid \phi (\omega ,\theta ) = r\}$ is at least
$2^{a},$ where $a$ is the cardinal of $\Nat$.}

\begin{proof} Define $H=\{h \mid  h \from \Nat\rightarrow R$ and
$h(0)=r\} $. By $(3),R$ has at least two distinct
elements, so the cardinal of $H$ is at least $2^{a}$.
Let $g \from \Nat\rightarrow\Omega$ be $1-1$ with
$g(0)=\omega $.
($g$ also exists by \eqn{c}.)
By \eqn{k}, for every $h \in H$ there exists $\theta_{h}$
such that 
$\forall n[\phi (g(n),\theta_{h}) = h(n)]$.
Note that, if $h\neq j,$ then
$\theta_{h}\neq \theta_{j},$ so the cardinal of
$\{\theta_{h} \mid h\in H\}$ is at least
$2^{a}$.
Setting $n=0$ yields 
$\forall h\in H[\phi (\omega ,\theta_{h})=r]$. \end{proof}

\lema f{If $\phi$ is rich, then for every 
$\omega\in \Omega$ there is an event $B\in \mathcal{B}$
satisfying $\omega\in B$ and $\pi (B) =0$.
Consequently, for every countable $C\subseteq \Omega$ there
is an event $B\in\mathcal{B}$ satisfying
$C\subseteq B$ and $\pi (B) =0$.}

\begin{proof} By countable additivity of $\pi$ it is sufficient to
establish, for arbitrary $\omega\in\Omega$ and for every
rational $x > 0,$ that there is an event $B\in\mathcal{B}$
satisfying $\omega\in\mathcal{B}$ and $\pi (B) < x$.
By \eqn{c}, there exist $r\in R$ and $X\in\mathcal{R}$ such
that $r\in X$ and $\mu (X) < 1$.
Let $f \from \Nat\rightarrow\Theta$ be $1-1$ with
$f(\Nat)\subseteq \{\theta \mid \phi (\omega ,\theta
)=r\}$.
Such a function exists by \lem{e}.
For $n\in\Nat,$ define $\alpha_{n} \from \Theta\rightarrow{\mathcal
R}$ by $\alpha_{n}(f(m))=X$ for all 
$m < n,$ and $\alpha_{n}(\theta )=R$ for every other 
$\theta $. Define $B_n = \Psi_\Omega^{-1}(\alpha_{n}^\#)$, and note
that $\omega \in B_n$.
By \eqn{h}, $\mu(B_n) = \kappa_{\pi}(\alpha_{n}^{\#})=(\mu (X))^{n}$.
For $n$ sufficiently large, 
$(\mu (X))^{n} < x$. \end{proof}

\theor B{If $\phi$ is rich, then $\b\Phi$
is a homogeneous sample-distribution limit for $\b{R}$.}

\begin{proof} $\b\Phi$ is a sample-distribution limit by \thm{A}, so only
\eqn{o} has to be verified.
That is, it must be shown that if $A\in\mathcal{F}$ and
$B\in \mathcal{R},$ then the set of $\omega$ satisfying
\display p{\nu (A\cap\phi_{\omega}^{-1}(B))=\nu
(A)\mu (B)}
\noindent is an event of $\mathcal{B}$ having probability 1.

That \eqn{p} holds a.s.~will first be proved in the case that 
$A\in \Psi_\Theta (\mathcal{A}_\mathcal{F})$.\footnote{$\Psi_\Theta$ is dual to $\Psi_\Omega$ defined in
  \eqn{l}. Specifically, first define $\mathcal{A}_{\mathcal{F}}^{0}$,
  $\mathcal{A}_{\mathcal{F}}$, and $\mathcal{E}_{\mathcal{F}}$ by
  dualizing equations \eqn{d}, \eqn{g}, and \eqn{C} respectively. Then
  define $\Psi_\Theta$ by dualizing \eqn{l}.}
By \lem{f}, it is sufficient to prove that \eqn{p} holds for all
but countably many $\omega$.
Moreover by (\eqnref{f}$'$), (\eqnref{g}$'$) and (\eqnref{l}$'$), it is sufficient to prove
\eqn{p} on this complement for
$A=\Psi_\Theta (\alpha^{\#}),$ where $\alpha\in\mathcal{A}_{\mathcal
F}^{0}$. 
That is, any element $A$ of $\Psi_\Theta (\mathcal{A}_\mathcal{F})$ is
a finite disjoint union of sets $\Psi_\Theta (\alpha^{\#})$.
Let $A=\Psi_\Theta (\alpha^{\#}),$ then, and define
$F=\{\omega \mid \alpha (\omega )\neq R\} $.
By (\eqnref{d}$'$), $F$ is finite.
Suppose that $\omega^{\ast}\not\in F,$ and define
$\beta \in \mathcal{A}_\mathcal{F}^{0}$ by 
$\beta (\omega^{\ast})=B$ and $\beta (\omega ) = R$
otherwise.
Also define $\gamma\in\mathcal{A}_\mathcal{F}^{0}$ by
$\gamma (\omega^{\ast})=B$ and
$\gamma (\omega ) = \alpha (\omega )$ otherwise.
Note that $\alpha^{\#}\cap\beta^{\#}=\gamma^{\#}$ that
$\kappa_{\nu}(\gamma^{\#})=\kappa_{\nu}(\alpha^{\#})\kappa_{\nu}(\beta^{\#})$
by (\eqnref{h}$'$), and that $\phi_{\omega^{\ast}}^{-1}(B)=\Psi_\Theta (\beta^{\#})$.
These facts establish \eqn{p} for $\alpha^{\#}$, $B,$ and all
$\omega^{\ast}\not\in F$.

Now using the fact that \eqn{p} holds a.s.~on
$\Psi_\Theta (\mathcal{A}_\mathcal{F}),$ it will be shown that \eqn{p}
holds a.s.~for arbitrary $A\in\mathcal{F}$.
Since $\mathcal{A}_\mathcal{F}$ generates $\mathcal{E}_\mathcal{F},$
lemma \rsltref{c}$'$ asserts that $\Psi_\Theta (\mathcal{A}_\mathcal{F})$ generates
$\mathcal{F}$.
Therefore, by Halmos ((1974), \S 13, Theorem D) and lemma \rsltref{a}$'$,
there exists a sequence
$\{ X_{n}\}_{n\in\Nat}\subseteq \mathcal{A}_\mathcal{F}$
such that
\display r{\forall n [\nu (\Psi_\Theta (X_{n})\bigtriangleup A) <
1/n].}
Let $B \in \mathcal{R}$, and define $Z_n = \{ \omega \mid \nu
(A\cap\phi_{\omega}^{-1}(B))\neq\nu (A)\mu (B) \}$. Since $X_{n}\in
\mathcal{A}_\mathcal{F}$, $\pi(Z_n)=0$. 

Replacing $A$ in \eqn{p} by $\Psi_\Theta (X_{n})$ and applying
\eqn{r} yields, for every $\omega \neq Z_n$ and $n>0$,
\display s{|\nu (A\cap\phi_{\omega}^{-1}(B)) - \nu (A)\mu (B)|
< 1/n}
which implies \eqn{p} for every $\omega \notin \bigcup_{n \in \Nat} Z_n$.
That is, for every $B \in \mathcal{R}$,  \eqn{p} holds a.s.~for $A$.
Since this is true for every $A \in \mathcal{F}$, \eqn{o} is satisfied. \end{proof}

\section{Existence of a rich function}

A rich function with range $R$ is now proved to exist for
suitably chosen sets $\Theta$ and $\Omega $.
This function is constructed by transfinite recursion, using
some basic results of cardinal arithmetic.\footnote{These topics are
covered, for instance, in Takeuti and Zaring (1982).
The facts about cardinal arithmetic that will be used are (a)
exact analogues of rules for manipulating sums, products and
exponents of natural numbers, (b) generalization of the
distributive law to transfinite addition, and (c) the facts that
the sum and the product of two infinite cardinals are both equal
to the maximum of the two operands, and that the cardinal of the
set of functions from one set to another is the cardinal of the
range taken to the exponent of the cardinal of the domain.}
$a$ and $r$ will denote the cardinals of $\Nat$ and $R$ respectively;
$\psi ,\rho ,\sigma$ and $\tau$ will denote ordinal
numbers; and $q$, $s$ and $t$ will denote cardinal
numbers.
Addition, multiplication and exponentiation will refer to
cardinal operations.

A rich function $\phi  \from \Omega\times\Theta\rightarrow R$ will
be obtained as the union of a transfinite nested sequence of
partial functions.
The sequence must be chosen to that $\phi$ will be a total
function that satisfies \eqn{k}.
At each stage of the sequence, either
$\phi (\omega ,\theta )$ will be defined for some specified
element $(\omega ,\theta )$ of the domain in order to
ensure that $\phi$ will be total, or else an instance of one
of the implications in \eqn{k} will be satisfied.
These characteristics of the function to be determined will be
called {\it features.}
The set $T$ of features is given by

\display t{T=\{ (h,f)\in R^\Nat\times\Theta^\Nat \mid f \text{ is
1--1}\}\cup\{ (h,g)\in R^\Nat\times\Omega^\Nat \mid g \text{ is
1--1}\}\cup[\Omega\times\Theta ].}

\noindent Let $N$ be an {\it enumeration} of $T$.
That is, suppose that $\tau$ is an ordinal and that,

\display u{N \from \tau\rightarrow T \text{ is onto } T.}

Let $P$ denote the set of partial functions from
$\Omega\times\Theta$ to $R$.
That is, $P$ is given by

\display q{P=\{ p \mid \exists A[A\subseteq \Omega\times
\Theta \text{ and } p \from A\rightarrow R]\}.}

\noindent Define the domain of a partial function, and the
projections of the domain on $\Theta$ and $\Omega$, by

\display v{ \begin{split}
  D(p) = \{ (\omega,\theta ) \mid \exists r \; p(\omega ,\theta) = r\} &\qquad
D_{\Omega}(p) = \{ \omega \mid \exists \theta(\omega ,\theta)\in D(p)\}\\
D_{\Theta}(p) = \{ \theta & \mid \exists\omega(\omega ,\theta)\in D(p)\}.
\end{split}}

A partial function {\it forces} a feature if it appropriately relates the
values of the two functions to which the feature refers.
To be precise,

\display w{p \text{ forces } (h,f)\in R^\Nat\times\Theta^\Nat
  \text{ iff } \exists\omega\forall n[p(\omega ,f(n)) = h(n)].}

\display x{p \text{ forces } (h,g)\in R^\Nat\times\Omega^\Nat
  \text{ iff } \exists\theta\forall n[p(g(n),\theta ) = h(n)].}

\display y{p \text{ forces } (\omega ,\theta )\in \Omega\times\Theta
 \text{ iff } \exists r[p(\omega ,\theta ) = r].}

A {\it fully specified sequence} of partial functions is a transfinite
nested sequence such that every feature is eventually forced. That
is, a function $S \from \tau\rightarrow P$ is a fully specified
sequence if

\display z{\forall\rho\forall\sigma[\rho <\sigma <\tau
 \implies  S(\rho )\subseteq S(\sigma )]} 

\noindent and

\display A{\forall \sigma < \tau[S(\sigma + 1)
\text{ forces } N(\sigma )].}

\lema g{If $S$ is a fully specified sequence
and $\phi = \bigcup_{\sigma <\tau}S(\sigma )$, then
$\phi$ is a rich function.}
 
\begin{proof} That $\phi$ is a total function satisfying \eqn{k} follows
directly from \eqn{t}--\eqn{A}. \end{proof}

It can be shown that if $S$ is a fully specified sequence, then the
domain of $\bigcup_{\sigma <\tau}S(\sigma)$ cannot have larger
cardinality than does $\Omega\times\Theta$.
Thus, the cardinality of $T$ implied by \eqn{t} must not be
larger than this. Otherwise the domain of $\phi$ would be exhausted
before all features had been forced.
The next lemma establishes that $\Theta$ and $\Omega$ can
be taken to be of a cardinality such that this problem will not
arise.

\lema h{Let $r$ be the cardinal of $R,$
and define $t=r^{a}$.\footnote{Note that, if $r=2^{a},$ then
$t=r$. $r=2^{a}$ if $(R, \mathcal{R})$ is a standard Borel space.}
Then $t^{a}=t$.
If $\Theta$ and $\Omega$ are of cardinality $t,$ then
$T$ is also of cardinality $t$.}

\begin{proof} The first assertion is true because $a^{2}=a,$ so that
$t^{a} =(r^{a})^{a}=r^{\left( a^{2} \right)}=r^{a}=t$.  To prove the
  second assertion, note first that $\{ (h,f)\in
  R^\Nat\times\Theta^{N} \mid f$ is $1-1\}$ is a subset of $T$.  This
  subset has at cardinality at least $t,$ since $t$ is the cardinal of
  $R^\Nat$.  $T$ is a subset of $(R^\Nat \times \Theta^\Nat) \cup
  (R^\Nat \times \Omega^\Nat)\cup (\Omega\times\Theta),$ the cardinal
  of which is $t \cdot t^{a} + t \cdot t^{a} +t \cdot t = t$
  also. \end{proof}

\lema i{Let $r$ be the cardinal of $R$, $t=r^{a}\!,$ $\Theta$ and $\Omega$ be
  disjoint sets of cardinality $t$, and $\tau$ be the initial ordinal
  of $t$.  Then there exist an enumeration $N \from \tau\rightarrow T$ and a
  fully specified sequence $S \from \tau\rightarrow P$.}

\begin{proof} The existence of $N$ follows from \lem{h}.
$S$ will now be described recursively.  That is, at each stage
  $\sigma$, the graph of $S(\sigma )\in P$ will be described.  Let
  $\sigma$ be the first ordinal for which $S$ has not been defined,
  and let $s$ be the cardinal of $\sigma$.  Consider the induction
  hypothesis that, if $\psi <\sigma$ and $q$ is the cardinal of
  $\psi$, then $D(S(\psi ))$ has cardinality no greater than $q+a,$ as
  well as that \eqn{z} and \eqn{A} are satisfied (with $\psi$ replacing
  $\tau )$.  If $\sigma =0$, define $S(\sigma) =\emptyset$.  At
  $0$, \eqn{z} and \eqn{A} are satisfied trivially.  If $\sigma =\rho +1$
  and $\rho$ satisfies the induction hypothesis, then $D(S(\rho
  )),D_{\Theta}(S(\rho ))$ and $D_{\Omega}(S(\rho ))$ all have
  cardinality no greater than $s+a < t$.  Therefore $\Theta\setminus
  D_{\Theta}(S(\rho ))$ and $\Omega\setminus D_{\Omega}(S(\rho ))$
  are nonempty.  Let $\theta^{\ast}\in\Theta\setminus
  D_{\Theta}(S(\rho )),$ $\omega^{\ast}\in \Omega\setminus
  D_{\Omega}(S(\rho )),$ and $r^{\ast}\in R$.  If $N(\rho )=(h,f)\in
  R^\Nat\times\Theta^\Nat,$ define $S(\sigma ) = S(\rho
  )\cup\{ (\omega^{\ast},f(n),h(n)) \mid n\in\Nat\} $.  If $N(\rho
  )=(h,g)\in R^\Nat\times \Omega^\Nat,$ define $S(\sigma ) =
  S(\rho )\cup\{ (g(n),\theta^{\ast},h(n)) \mid n\in\Nat\} $.  If
  $N(\rho ) = (\omega ,\theta )\in\Omega\times\Theta ,$ define
  $S(\sigma ) = S(\rho )$ if $(\omega ,\theta )\in D(S(\rho ))$ and
  define $S(\sigma ) = S(\rho )\cup\{ (\omega ,\theta ,r^{\ast})\}$
  otherwise.  For each of the three types of feature, $S(\rho )$ is
  extended to a countable set outside $D(S(\rho ))$ to obtain
  $S(\sigma )$.  Thus \eqn{z} and \eqn{A} hold, and $D(S(\sigma ))$ has
  cardinality no greater than $s+a$.  If $S$ has been defined up to
  $\sigma$ and $\sigma$ is a limit ordinal, then define $S(\sigma)
  =\bigcup_{\rho<\sigma}S(\rho )$.  Equations \eqn{z} and \eqn{A} continue to
  hold.  $S(\sigma )$ is a union of $s$ sets of cardinality no greater
  than $s+a,$ so its cardinality is no greater than $s\cdot(s+a)=s+a$.
  By transfinite induction, then, $S \from \tau\rightarrow P$ is defined and
  satisfies \eqn{z} and \eqn{A}. \end{proof}

\theor C{A rich function exists, and $\b{R}$
possesses a homogeneous sample-\hfill\break distribution limit.}

\begin{proof} The existence of a rich function follows from \lem{g} and
\lem{i}.
Given this function, the existence of a homogeneous sample-distribution limit
for $\b{R}$ follows from \thm{B}. \end{proof}

\section{Nonatomicity and nonmeasurability of homogeneous
sample-distribution limits}

The construction of $\phi$ has guaranteed that all of the
sections of $\phi $, both with respect to $\Theta$ and to
$\Omega ,$ are measurable.
However, the measurability of $\phi$ with respect to the
product $\sigma$-algebra $\mathcal{B}\times\mathcal{F}$ has not
been asserted.
In this section it will be shown that, under a mild restriction,
no homogeneous sample-distribution limit can be jointly
measurable in its two variables.
The restriction in that $\mathcal{R}$ should contain a set of
$\mu$-measure strictly between 0 and 1.\footnote{If $\mathcal{R}$
is contained in the completion by measure-zero sets of the
$\sigma$-field of invariant sets of an ergodic transformation
on $R,$ then (by definition) the restriction is not
satisfied.
This seems to be the only nontrivial case of practical interest
in which the restriction would not be satisfied.
Doob (1937) has proved the nonmeasurability of $\phi$ when
$\b\Theta$ is taken to be the unit interval with Lebesgue
measure, but his proof does not generalize to the situation where
$\b\Theta$ is not countably generated.}

The proof of this result makes use of the fact that, under the
restriction, every homogeneous sample-distribution limit has a
nonatomic population measure (in the sense of Halmos ((1974), \S
40): that every set of positive measure has a subset of strictly
smaller positive measure).
This fact, which is of some independent interest, is now proved.
Note that lemma \rsltref{f}$'$ fails to imply that the population measure
is nonatomic because the $\sigma$-algebra on which it is
defined is not countably generated. (cf. the example at the end
of this section.)

\lema j{If $\phi$ satisfies $(16),$ and if
$\mathcal{R}$ contains a set of $\mu$-measure strictly between
0 and 1, then $\nu$ is nonatomic.}

\begin{proof} Suppose that $B\in\mathcal{R}$ and
$0 < \mu (B) < 1$.
  Then \eqn{o} implies that if $0 < \nu (A) < 1$, then $0 <
  \nu(A \cap \phi_\omega^{-1}(B)) < \nu (A)$. Thus $\nu$ cannot
  have an atom. \end{proof} 

\theor D{Suppose that $\phi$ satisfies
$(1)$ and satisfies $(16),$ and that $\mathcal{R}$ contains a
set of $\mu$-measure strictly between 0 and 1.
Then $\phi$ is not measurable with respect to 
$\mathcal{B}\times\mathcal{F}$.}

\begin{proof} Suppose that $A\in\mathcal{R}$, $\mu (A) = a$, and
$0<a<1$.
It will be assumed that $\phi^{-1}(A)$ is measurable, and
this assumption will be shown to lead to a contradiction.
Define $\psi  =\pi\times\nu$ and 
$b=(a-a^{2})/2$.
Since $0 < a < 1$, $b > 0$.
If $\phi^{-1}(A)$ is measurable, then there exist a finite
set $I$ and sets $\{ B_{i}\}_{i\ \in I}\subseteq{\mathcal
B}$ and $\{ F_{i}\}_{i\in I}\subseteq \mathcal{F}$ such
that  $i\neq j \implies [(B_{i}\times
F_{i})\cap(B_{j}\times F_{j}) = \emptyset ]$ and
$\psi (G) < b,$ where\hfill\break
$G=\phi^{-1}(A)\bigtriangleup\bigcup_{i \in I}(B_{i}\times
F_{i})$.
(Halmos (1974), \S 33, Theorem E  and \S 13, Theorem D).

By Fubini's theorem (Halmos (1974), \S 36, Theorem B)
$\psi (G)=\int_{\Theta}\pi
(\phi_{\theta}^{-1}(A)\bigtriangleup\bigcup \{B_{i} \mid \theta\in
F_{i}\} )d\nu (\theta )$.
Therefore, for some $H\in\mathcal{F}$, $\nu (H) > 0$
and $\forall\theta\in H\; \pi
(\phi_{\theta}^{-1}(A)\bigtriangleup\bigcup \{ B_{i} \mid \theta\in
F_{i}\} ) < b$.
Define an equivalence relation $\approx$ on
$H$ by $\eta\approx\theta$ iff 
$\forall i\in I\; \eta\in F_{i}\equiv \theta\in F_{i}$.
This relation induces a finite measurable partition of $H,$
and $H$ is infinite because $\nu$ is nonatomic by \lem{j}.
Thus there exist distinct $\eta$ and $\theta$ such that
$\eta\approx\theta $.
Define 
$J=\bigcup\{ B_{i} \mid \eta\in F_{i}\} =\bigcup \{
B_{i} \mid \theta\in F_{i}\} $.
Then
\medskip
\item{\eqn{A}} $\pi (\phi_{\eta}^{-1}(A)\bigtriangleup J) < b$ and
$\pi (\phi_{\theta}^{-1}(A)\bigtriangleup J) < b$.
\medskip
\noindent Since
$\pi (\phi_{\eta}^{-1}(A)) = \pi (\phi_{\theta}^{-1}(A)) = a$
by \eqn{a}, \eqn{A} implies that

\display B{\pi
  (J\cap\phi_{\eta}^{-1}(A)\cap\phi_{\theta}^{-1}(A)) 
> a-2b = a^{2}.}

\noindent This contradicts the independence of $\phi_{\eta}$
and $\phi_{\theta},$ which requires that\\
$\pi (\phi_{\eta}^{-1}(A)\cap\phi_{\theta}^{-1}(A)) =
a^{2}$\end{proof}

The hypothesis of \thm{D} concerning $\mathcal{R}$ (ensuring
that $\mu$ is nonatomic) evidently cannot be dropped.
That is, if $\mu (\{ r\} )=1$ and
$\forall\omega\,\forall\theta\:\phi (\omega ,\theta )=r,$
then $\phi$ is measurable and satisfies $(1)$ and \eqn{o}.
A further example will show that the hypothesis cannot be
weakened to the statement that $\mu$ is not concentrated at a
single point.
Let $\mathcal{E}$ be the $\sigma$-algebra of countable and
co-countable subsets of $[0,1],$ and define
$\varepsilon (A)$ to be 0 if $A$ is countable and 1 if 
$X\setminus A$ is countable.
Let $\b\Omega = \b\Theta = \b{R} = ([0,1],\mathcal{E},\varepsilon
)$.
Then any 1--1 function
$\phi  \from [0,1]^{2}\rightarrow [0,1]$ is product measurable
and satisfies \eqn{a} and \eqn{o}.

\bigskip

\parindent=0pt
\centerline{\bf References}

\hangindent=36pt Anderson, R.~M., (1991), ``Nonstandard Analysis with 
Applications to Economics,'' in {\it Handbook of Mathematical Economics},
Volume IV. K.~J.~Arrow and M.~D.~Intrilligator, eds., North Holland, 
Amsterdam.

\hangindent=36pt Doob, J.~L., (1937), ``Stochastic Processes
Depending on a Continuous Parameter,'' {\it Trans. Am. Math.
Soc.}, {\bf 42}, 107--140.

\hangindent=36pt Doob, J.~L., (1947), ``Probability in Function
Space,'' {\it Bull. Am. Math. Soc.}, {\bf 53}, 15--30.

\hangindent=36pt Doob, J.~L., (1953), {\it Stochastic Processes}, Wiley,
New York.

\hangindent=36pt Feldman, M., and C. Gilles, (1985), ``An
Expository Note on Individual Risk Without Aggregate Uncertainty,''
{\it J. Econ. Theory}, {\bf 35}, 26--322.

\hangindent=36pt Green, E.~J., (1987), ``Lending and the Smoothing of
Uninsurable Income,'' in {\it Intertemporal Trade and Financial
Intermediation}, E. Prescott and N. Wallace, eds., University of
Minnesota Press, Minneapolis.

\hangindent=36pt Halmos, P., (1974), {\it Measure Theory}, Springer, New
York.

\hangindent=36pt Judd, K. L., (1985), ``The Law of Large Numbers
with a Continuum of i.i.d.~Random Variables,'' {\it J. Econ.
Theory}, {\bf 35}, 19--25.

\hangindent=36pt Keisler, H.~J., (1977), ``Hyperfinite Model Theory,''
in {\it Logic Colloquium 1976}, R.~O.~Gandy and J.~M.~E.~Hylland,
eds., North Holland, Amsterdam.

\hangindent=36pt Parthasarathy, K.R., (1967), {\it Probability Measures
on Metric Spaces}, Academic Press, New York.

\hangindent=36pt Takeuti, G., and W. Zaring, (1982), {\it Introduction to
Axiomatic Set Theory}, 2nd ed., Springer, New York.

\hangindent=36pt Uhlig, H., (1996), ``A Law of Large Numbers for
Large Economies,'' {\it Econ. Theory}, {\bf 8}, 41--50.

\end{document}